\documentclass[12pt]{article}

\usepackage{amsmath,amsthm,amsfonts,latexsym,amsopn,verbatim,amscd,amssymb,}
\usepackage{amssymb}
\usepackage{amsfonts}

\usepackage{graphics,color}
\usepackage{graphicx}
\usepackage{amssymb,array,arydshln}
\usepackage{amsfonts}
\usepackage{amsbsy}
\usepackage{amssymb}
\usepackage{amsmath}
\usepackage{cite}
\usepackage{pst-all}
\usepackage{pstricks}
\usepackage{graphics}
\usepackage{float}

\theoremstyle{plain}
\newtheorem{Thm}{Theorem}[section]
\newtheorem{Lem}[Thm]{Lemma}

\newtheorem{Cor}[Thm]{Corollary}
\theoremstyle{definition}

\newtheorem{Def}[Thm]{Definition}

\numberwithin{equation}{section}

\DeclareMathOperator{\idem}{Idem}
\DeclareMathOperator{\nil}{Nil}

\begin{document}
\begin{center}
\textbf{NIL CLEAN GRAPH OF RINGS }\\
\end{center}
\begin{center}
Dhiren Kumar Basnet\\
\small{\it Department of Mathematical Sciences, Tezpur University,
 \\ Napaam, Tezpur-784028, Assam, India.\\
Email: dbasnet@tezu.ernet.in}
\end{center}
\begin{center}
Jayanta Bhattacharyya \\
\small{\it Department of Mathematical Sciences, Tezpur University,
 \\ Napaam, Tezpur-784028, Assam, India.\\
Email: jbhatta@tezu.ernet.in}
\end{center}
\noindent \textit{\small{\textbf{Abstract:}  }} In this article, we have defined nil clean graph of a ring $R$. The vertex set is the ring $R$, two ring elements $a$ and $b$ are adjacent if and only if $a + b$ is nil clean in $R$. Graph theoretic properties like girth, dominating set, diameter etc. of nil clean graph have been studied for  finite commutative rings.
\bigskip

\smallskip

\smallskip

\bigskip
\section{Introduction}
In past, lot of mathematicians have studied graph theoretic properties of graph associated with rings. By graph, in this article we mean simple undirected graph. If $G$ denotes a graph, let $V(G)$ denotes the set of vertices and $E(G)$ denotes the set of edges. In $1988$, Istavan Beck \cite{cocr}, studied colouring of graph of a finite commutative ring $R$, where vertices are elements of $R$ and $xy$ is an edge if and only if $xy = 0$. For a positive integer $n$, let $\mathbb{Z}_n$ be the ring of integers modulo $n$. P. Grimaldi\cite{gfr} defined and studied various properties of the unit graph $G(\mathbb{Z}_n)$, with vertex set $\mathbb{Z}_n$ and two distinct vertices $x$ and $y$ are adjacent if $x+ y$ is a unit. Further in \cite{ug}, authors generalized $G(\mathbb{Z}_n)$ to unit graph $G(R)$, where $R$ is an arbitrary associative ring with non zero identity. They have studied properties like diameter, girth chromatic index etc.

In this paper we have introduced nil clean graph $G_N(R)$ associated with a finite commutative ring $R$. The properties like girth, diameter, dominating sets etc. of $G_N(R)$ have been studied.

Ring $R$ in this article, is finite commutative ring with non zero identity. An element $r$ of ring $ R$ is said to be \textit{nil clean element}\cite{nc, ajd} if for an idempotent $e$ in $R$ and a nilpotent $n$ in $R$, $r = n + e$. The ring $R$ is said to be \textit{nil clean} if each element of $R$ is nil clean. The set of nil clean elements of $R$ is denoted by $NC(R)$. The sets of idempotents and nilpotents of $R$  are denoted by $\idem(R)$ and $\nil(R)$ respectively. A ring $R$ is said to be a \textit{weak nil clean ring}\cite{cwnc}, if for each $r \in R$, $r = n + e$ or $r = n - e$, for some $n\in\nil(R)$  and $e\in \idem(R)$.
\section{Basic properties}
In this section we will define nil clean graph of a finite commutative ring and discuss its basic properties.

\begin{Def}
The \textit{nil clean graph} of a ring $R$, denoted by $G_N(R)$, is defined by setting $R$ as vertex set and defining two distinct verities $x$ and $y$ to be adjacent if and only if $x + y$ is a nil clean element in $R$. Here we are not considering loops at a point (vertex) in the graph.
\end{Def}
For illustration following is the nil clean graph of $GF(25)$, where $GF(25)$ is the finite field with $25$ elements.
\begin{align*}
GF(25) &\cong \mathbb{Z}_5[x] /\langle x^2 + x + 1\rangle\\
         & = \{ ax + b + \langle x^2 + x + 1\rangle : a, b \in \mathbb{Z}_5\}
\end{align*}
let us define $\alpha := x + \langle x^2 + x + 1\rangle$, then we have
$GF(25) = \{ \overline 0, \overline 1, \overline 2, \overline 3, \overline 4, \alpha, 2\alpha, 3\alpha, 4\alpha, 1 + \alpha, 1 + 2\alpha, 1 + 3\alpha, 1 + 4\alpha,
2 + \alpha, 2 + 2\alpha, 2 + 3\alpha, 2 + 4\alpha, 3 + \alpha, 3 + 2\alpha, 3 + 3\alpha, 3 + 4\alpha, 4 + \alpha, 4 + 2\alpha, 4 + 3\alpha, 4 + 4\alpha\}$. Observe that $NC(GF(25)) =\{\overline 0, \overline 1\}$.

\begin{figure}[H]  
\begin{pspicture}(0,4)(0,-4.5)
\scalebox{1}{
\rput(2,0){
\psdot[linewidth=.05](1,3)
\psdot[linewidth=.05](3,3)
\psdot[linewidth=.05](5,3)
\psdot[linewidth=.05](7,3)
\psdot[linewidth=.05](9,3)
\psline(1,3)(9,3)
\rput(1,3.3){$\overline 0$}
\rput(3,3.3){$\overline 1$}
\rput(5,3.3){$\overline 4$}
\rput(7,3.3){$\overline 2$}
\rput(9,3.3){$\overline 3$}

\psdot[linewidth=.05](1,2)
\psdot[linewidth=.05](3,2)
\psdot[linewidth=.05](5,2)
\psdot[linewidth=.05](7,2)
\psdot[linewidth=.05](9,2)
\psdot[linewidth=.05](1,0)
\psdot[linewidth=.05](3,0)
\psdot[linewidth=.05](5,0)
\psdot[linewidth=.05](7,0)
\psdot[linewidth=.05](9,0)
\rput(1,2.3){$\alpha$}
\rput(3,2.3){$4\alpha + 1$}
\rput(5,2.3){$\alpha + 4$}
\rput(7,2.3){$4\alpha + 2$}
\rput(9,2.3){$\alpha + 3$}
\rput(1,-.3){$4\alpha$}
\rput(3,-.3){$\alpha + 1$}
\rput(5,-.3){$4\alpha + 4$}
\rput(7,-.3){$\alpha + 2$}
\rput(9,-.3){$4\alpha + 3$}
\psline(1,2)(1,0)
\psline(9,0)(9,2)
\psline(1,0)(9,0)
\psline(1,2)(9,2)

\psdot[linewidth=.05](1,-1.5)
\psdot[linewidth=.05](3,-1.5)
\psdot[linewidth=.05](5,-1.5)
\psdot[linewidth=.05](7,-1.5)
\psdot[linewidth=.05](9,-1.5)
\psdot[linewidth=.05](1,-3.5)
\psdot[linewidth=.05](3,-3.5)
\psdot[linewidth=.05](5,-3.5)
\psdot[linewidth=.05](7,-3.5)
\psdot[linewidth=.05](9,-3.5)
\rput(1,-1.2){$2\alpha$}
\rput(3,-1.2){$3\alpha + 1$}
\rput(5,-1.2){$2\alpha + 4$}
\rput(7,-1.2){$3\alpha + 2$}
\rput(9,-1.2){$2\alpha + 3$}
\rput(1,-3.8){$3\alpha$}
\rput(3,-3.8){$2\alpha + 1$}
\rput(5,-3.8){$3\alpha + 4$}
\rput(7,-3.8){$2\alpha + 2$}
\rput(9,-3.8){$3\alpha + 3$}
\psline(1,-1.5)(1,-3.5)
\psline(9,-1.5)(9,-3.5)
\psline(1,-1.5)(9,-1.5)
\psline(1,-3.5)(9,-3.5)

}}
\end{pspicture}
\caption{Nil clean graph of $GF(25)$}\label{gf25}
\end{figure}
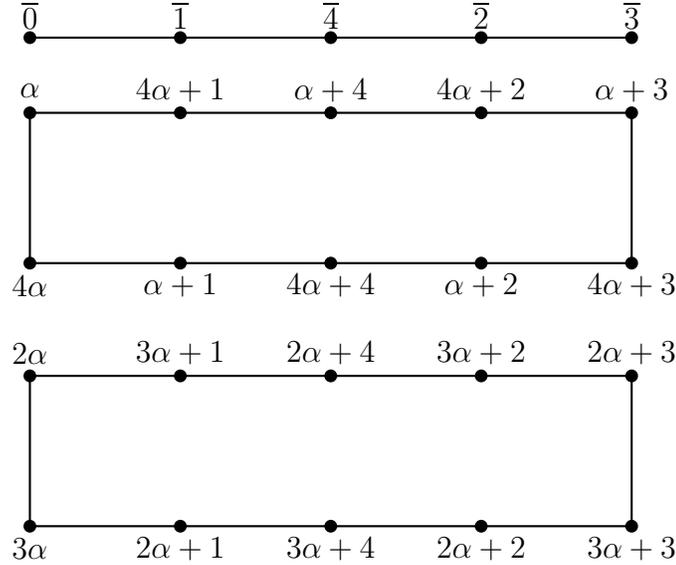

In graph theory, a \textit{complete graph} is a simple undirected graph in which every pair of distinct vertices is connected by a unique edge. So by this definition the following theorem follows.

 \begin{Thm}
 The nil clean graph $N_G(R)$ is a complete graph if and only if $R$ is a nil clean ring.
\end{Thm}
\noindent$Proof.$ Let $N_G(R)$ be complete nil clean graph of a ring $R$. For $r \in R$, $r$ is adjacent to $0$ implies $r = r + 0$ is nil clean, hence $R$ is nil clean. Converse is clear from definition of nil clean graph.$\hfill \square$

Two graphs $G_1$ and $G_2$ are said to be isomorphic if there exists an isomorphism from $G_1$ to $G_2$ we mean a bijective mapping $f:V(G_1)\rightarrow V(G_2)$ from $V(G_1)$ onto $V(G_2)$ such that two vertices $u_1$ and $v_1$ are adjacent in $G_1$ if and only if the vertices φ$f(u_1)$ and φ$f(v_1)$ are adjacent in $G_2$ \cite{igt}. For rings $R$ and $S$ if $R \cong S$ it is easy to see that $G_N(R) \cong G_N(S)$.

\begin{Lem}
Let $R$ be a commutative ring and idempotents lift modulo $\nil(R)$. If $x + \nil(R)$ and $y + \nil(R)$ are adjacent in $G_N(R/\nil(R))$ then every element of $x + \nil(R)$ is adjacent to every element of $y + \nil(R)$ in the nil clean graph $G_N(R)$.
\end{Lem}

\noindent$Proof.$ Let $x + \nil(R)$ and $y + \nil(R)$ be adjacent in $G_N(R/\nil(R))$, then
$$(x + \nil(R)) + (y + \nil(R)) = e + \nil(R),$$
where $e$ is an idempotent in $R$, as idempotents lift modulo $\nil(R)$. Thus from above we have
$x + y = e + n$, for some $n \in \nil(R)$ and hence $x$ and $y$ are adjacent in $G_N(R)$. Now for $a \in x + \nil(R)$ and $b \in y + \nil(R)$, we have $a = x + n_1$ and $b = y + n_2$, for some $n_1$, $n_2 \in \nil(R)$. Therefore $a + b = e + (n - n_1 - n_2)$. Hence, $a$ and $b$ are adjacent in $G_N(R)$.$\hfill \square$

Let $G$ be a graph, for $x\in V(G)$, the degree of $x$ denoted by $deg(x)$, is defined to be number of edges of $G$ whose one of the end point is $x$. The neighbour set of $x \in V(G)$, is denoted and defined as by $N_G(x):=\{ y\in V(G) | y\mbox{ is adjacent to } x\}$ and the set $N_G[x]= N_G(x) \cup \{x\}$.

\begin{Lem}\label{deg}
Let $G_N(R)$ be the nil clean graph of a ring $R$. For $x \in R$ we have the following:
\begin{enumerate}
  \item If $2x$ is nil clean, then $deg(x) = |NC(R)| - 1$.
 \item If $2x$ is not nil clean, then $deg(x) = |NC(R)|$.
\end{enumerate}
\end{Lem}

\noindent$Proof.$ Let $x \in R$, observe that $x + R = R$, so for every $y \in NC(R)$, there exists a unique element $x_y \in R$, such that $x + x_y = y$. Thus we have $deg(x) \leq |NC(R)|$. Now if $2x \in NC(R)$, define $f : NC(R) \rightarrow N_{G_N(R)}[x]$ given by $f(y) = x_y$. It is easy to see that $f$ is a bijection and therefore $deg(x) = |N_{G_N(R)}(x)| = |N_{G_N(R)}[x]| - 1 = |NC(R)| -1$.
If $2x \notin  NC(R)$, define $f : NC(R) \rightarrow N_{G_N(R)}(x)$ given by $f(y) = x_y$. Then $f$ is a bijection and therefore $deg(x) = |N_{G_N(R)}(x)| =  |NC(R)|$.$\hfill \square$

A graph $G$ is said to be \textit{connected} if for any two distinct vertices of $G$, there is a path in $G$ connecting them.

\begin{Thm}
 For a ring $R$. The following hold :
 \begin{enumerate}
   \item $G_N(R)$ need not be connected.
   \item Let $R = \mathbb{Z}_n$. For $\overline{a} \in \mathbb{Z}_n$ there is a path from $\overline{a}$ to $\overline{0}$.
   \item $G_N(\mathbb{Z}_n)$ is connected.
   \item Let $R = \mathbb{Z}_n$. For $A \in M_n(\mathbb{Z}_n)$ there is a path from $A$ to $0$, where $0$ is the zero matrix of $M_n(\mathbb{Z}_n)$.
   \item $G_N(M_n(\mathbb{Z}_n))$ is connected.
 \end{enumerate}
 \end{Thm}

\noindent$Proof.$ $(i)$ Clear by the graph $G_N(GF(25))$, figure \ref{gf25}. For $(ii)$ and $(iii)$ replacing $p$ by $n$ in the figure \ref{zp}, we get a Hamiltonian path in $G_N(\mathbb{Z}_n)$. Now for the proof of $(iv)$, let $ A = [a_{ij}] \in M_n(Z_n)$, now we define $A_1= [a1_{ij}]$ where $a1_{ij} = -a_{ij}$ for $i \geq j$, otherwise $a1_{ij} =0$, observe that $A_1 + A$ is nilpotent hence nil clean. Thus there exists an edge between $A$ and $A_1$. Again define $A_2 = [a2_{ij}]$ where $a2_{ij} = a_{ij}$ for $i = j$, otherwise $a1_{ij} =0$, then clearly we have an edge between $A_1$ and $A_2$ in $G_N(M_n(\mathbb{Z}_n))$. For each element $a_{ii}$ of $A_2$, by $(ii)$ we have a path
$\{ a_{ii}, b_{i1}, b_{i2}, b_{i3},\dots, b_{ik_i} = \overline 0\}$ of length $k_i \in \mathbb{N}$ to $\overline 0$. Now let $K = \max\{k_i: 1\leq i\}$, we will construct a path of length $K$ from $A_2$ to $0$, as follows. Define $B_{i}=[b1_{jl}]$, where $b1_{jj}=b_{ji}$, if $b1_{ij}$ appears in some above paths, otherwise $b1_{ij}=0$, $1\leq i \leq K$. Thus $\{ A, A_1, A_2, B_1, B_2, \dots, B_K = 0\}$ is a path from $A$ to $0$ in $G_N(M_n(\mathbb{Z}_n))$. Lastly $(v)$ follows from $(iv)$.$\hfill \square$

Following result related to finite field will be used in upcoming results, we state it here to maintain the self containment of the article.

\begin{Lem}\label{F}
A ring $R$ is a finite commutative reduced ring with no non trivial idempotents if and only if $R$ is a finite field.
\end{Lem}

\noindent$Proof. (\Rightarrow)$ Let $0 \neq x \in R$. Observe the set $A = \{ x^k: k\in \mathbb{N}\}$ is a finite set. Therefore there exist $m > l$ such that $x^l = x^m$. Note that
$ x^l  = x^m = x^{m-l+l} = x^{m-l}.x^l= x^{m-l}.x^m = x^{2m-l+l-l} = x^{2(m-l)+l} =\dots = x^{k(m-l)+l}$, where $k$ is a natural number. Now we have $[x^{l(m-l)}]^2 = x^{l(m-l)}.x^{l(m-l)} = x^{l(m-l) + l(m-l) + l - l}\ = x^{l(m-l)+l}.x^{l(m-l)-l}= x^l.x^{l(m-l)-l} = x^{l(m-l)} $
that is $x^{l(m-l)}$ is an idempotent. Thus $x^{l(m-l)} = 1$ which gives that $x$ is a unit, therefore $R$ is a finite field. $(\Leftarrow)$ Obvious.$\hfill \square$

\section{Invariants of nil clean graph }
In this section, we prove some results related to invariants of graph theory. Following subsection is for girth of $G_N(R)$.
\subsection{Girth of $G_N(R)$}
For a graph $G$, the \textit{girth} of $G$ is the length of the shortest cycle in $G$. We have following results on girth of nil clean graph.

\begin{Thm}\label{girth}
The following hold for nil clean graph $G_N(R)$ of $R$:
 \begin{enumerate}
  \item If $R$ is not a field, then girth of $G_N(R)$ is equal to 3.
  \item  If $R$ is a field
   \begin{enumerate}
   \item[(a)] girth is $2p$ if $R \cong GF(p^k)$ (field of order $p^k$), where $p$ is a odd prime and $k > 1$;
   \item[(b)] girth is infinite, in fact $G_N(R)$ is a path, otherwise.
   \end{enumerate}
\end{enumerate}
\end{Thm}
\noindent$Proof.$ (i) Let $R$ has at least one non-trivial idempotent or non trivial nilpotent. If $e \in R$ be an nontrivial idempotent, then we have

\begin{figure}[H]
\begin{pspicture}(0,-1)(0,2.2)
\centering
\scalebox{1}{
\rput(-2,0){
\psdot[linewidth=.05](7,1)
\psdot[linewidth=.05](9,2)
\psdot[linewidth=.05](9,0)
\psline(7,1)(9,2)(9,0)(7,1)
\rput(6.7,1){$0$}
\rput(9.7,2){$(1-e)$}
\rput(9.4,0){$e$}
}}
\end{pspicture}
\caption{}
\end{figure}

so the girth of $G_N(R)$ is 3. Again if $R$ contains a nontrivial nilpotent $n \in R$, then we have the cycle

\begin{figure}[H]
\begin{pspicture}(0,-1)(0,2.2)
\centering
\scalebox{1}{
\rput(-2,0){
\psdot[linewidth=.05](7,1)
\psdot[linewidth=.05](9,2)
\psdot[linewidth=.05](9,0)
\psline(7,1)(9,2)(9,0)(7,1)
\rput(6.7,1){$0$}
\rput(9.7,2){$1$}
\rput(9.4,0){$n$}
}}
\end{pspicture}
\caption{}
\end{figure}

so the girth is $3$. By \textbf{Lemma \ref{F}} rings without non trivial idempotents and nilpotents are field. This proves (i).

(ii) As set of nil clean elements of a finite field is $\{ 0, 1\}$, so the nil clean graph of
$\mathbb{F}_p$, where $p$ is a prime, is
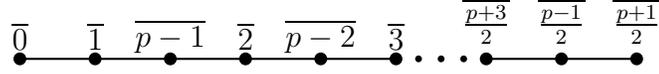
\begin{figure}[H]
\begin{pspicture}(0,0)(0,1)
\scalebox{1}{
\rput(2,0){
\psdot[linewidth=.05](1,0)
\psdot[linewidth=.05](2,0)
\psdot[linewidth=.05](3,0)
\psdot[linewidth=.05](4,0)
\psdot[linewidth=.05](5,0)
\psdot[linewidth=.05](6,0)
\psdot[linewidth=.01](6.3,0)
\psdot[linewidth=.01](6.6,0)
\psdot[linewidth=.01](6.9,0)
\psdot[linewidth=.05](7.2,0)
\psdot[linewidth=.05](8.2,0)
\psdot[linewidth=.05](9.2,0)
\psline(1,0)(2,0)(3,0)(4,0)(5,0)(6,0)
\psline(7.2,0)(8.2,0)(9.2,0)
\rput(1,0.3){$\small{\overline 0}$}
\rput(2,0.3){$\small{\overline 1}$}
\rput(3,0.3){$\small{\overline{p-1}}$}
\rput(4,0.3){$\small{\overline 2}$}
\rput(5,0.3){$\small{\overline{p-2}}$}
\rput(6,0.3){$\small{\overline 3}$}
\rput(7.2,0.5){$\small{\overline{\frac{p+3}{2}}}$}
\rput(8.2,0.5){$\small{\overline {\frac{p-1}{2}}}$}
\rput(9.2,0.5){$\small{\overline{\frac{p+1}{2}}}$}
}}
\end{pspicture}
\caption{Nil clean graph of $\mathbb{Z}_p$}\label{zp}
\end{figure}

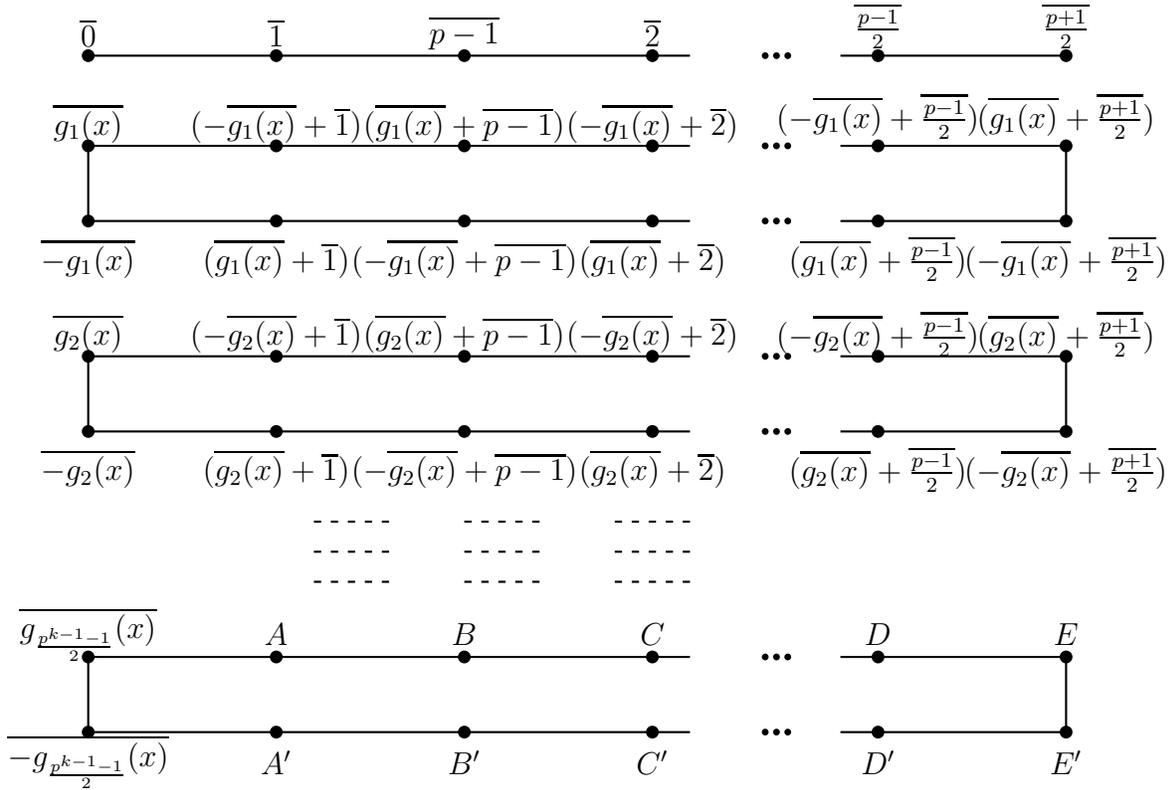
\begin{figure}[H]  
\begin{pspicture}(0,-1.5)(0,8.5)
\scalebox{1}{
\rput(0,0){
\rput(0,8){
\psdot[linewidth=.05](0,0)
\psdot[linewidth=.05](2.5,0)
\psdot[linewidth=.05](5,0)
\psdot[linewidth=.05](7.5,0)
\psdot[linewidth=.01](9,0)
\psdot[linewidth=.01](9.15,0)
\psdot[linewidth=.01](9.3,0)
\psdot[linewidth=.05](10.5,0)
\psdot[linewidth=.05](13,0)
\rput(0,.3){
\rput(0,0){$\small{\overline 0}$}
\rput(2.5,0){$\small{\overline 1}$}
\rput(5,0){$\overline{p-1}$}
\rput(7.5,0){$\small{\overline 2}$}
\rput(0,.1){
\rput(10.5,0){$\small{\overline {\frac{p-1}{2}}}$}
\rput(13,0){$\small{\overline{\frac{p+1}{2}}}$}}}
\psline(0,0)(8,0)
\psline(10,0)(13,0)
}

\rput(0,6.8){
\psdot[linewidth=.05](0,0)
\psdot[linewidth=.05](2.5,0)
\psdot[linewidth=.05](5,0)
\psdot[linewidth=.05](7.5,0)
\psdot[linewidth=.01](9,0)
\psdot[linewidth=.01](9.15,0)
\psdot[linewidth=.01](9.3,0)
\psdot[linewidth=.05](10.5,0)
\psdot[linewidth=.05](13,0)
\rput(0,.3){
\rput(0,0){$\small{\overline{g_1(x)}}$}
\rput(2.5,0){$\small{(-\overline{g_1(x)}+\overline 1})$}
\rput(5,0){$\small{(\overline{g_1(x)}+\overline{p-1})}$}
\rput(7.5,0){$\small{(-\overline{g_1(x)}+\small{\overline 2})}$}
\rput(10.5,0.1){$\small{(-\overline{g_1(x)}+\overline {\frac{p-1}{2}})}$}
\rput(13,0.1){$\small{(\overline{g_1(x)}+\overline{\frac{p+1}{2}})}$}}
\psline(0,0)(8,0)
\psline(10,0)(13,0)
\rput(0,-1){
\psdot[linewidth=.05](0,0)
\psdot[linewidth=.05](2.5,0)
\psdot[linewidth=.05](5,0)
\psdot[linewidth=.05](7.5,0)
\psdot[linewidth=.01](9,0)
\psdot[linewidth=.01](9.15,0)
\psdot[linewidth=.01](9.3,0)
\psdot[linewidth=.05](10.5,0)
\psdot[linewidth=.05](13,0)
\rput(0,-.5){
\rput(0,0){$\small{\overline{-g_1(x)}}$}
\rput(2.5,0){$\small{(\overline{g_1(x)}+\overline 1})$}
\rput(5,0){$\small{(-\overline{g_1(x)}+\overline{p-1})}$}
\rput(7.5,0){$\small{(\overline{g_1(x)}+\small{\overline 2})}$}
\rput(10.5,0){$\small{(\overline{g_1(x)}+\overline {\frac{p-1}{2}})}$}
\rput(13,0){$\small{(-\overline{g_1(x)}+\overline{\frac{p+1}{2}})}$}}
\psline(0,0)(8,0)
\psline(10,0)(13,0)}
\psline(0,0)(0,-1)
\psline(13,0)(13,-1)
}
\rput(0,4){
\psdot[linewidth=.05](0,0)
\psdot[linewidth=.05](2.5,0)
\psdot[linewidth=.05](5,0)
\psdot[linewidth=.05](7.5,0)
\psdot[linewidth=.01](9,0)
\psdot[linewidth=.01](9.15,0)
\psdot[linewidth=.01](9.3,0)
\psdot[linewidth=.05](10.5,0)
\psdot[linewidth=.05](13,0)
\rput(0,.3){
\rput(0,0){$\small{\overline{g_2(x)}}$}
\rput(2.5,0){$\small{(-\overline{g_2(x)}+\overline 1})$}
\rput(5,0){$\small{(\overline{g_2(x)}+\overline{p-1})}$}
\rput(7.5,0){$\small{(-\overline{g_2(x)}+\small{\overline 2})}$}
\rput(10.5,0){$\small{(-\overline{g_2(x)}+\overline {\frac{p-1}{2}})}$}
\rput(13,0){$\small{(\overline{g_2(x)}+\overline{\frac{p+1}{2}})}$}}
\psline(0,0)(8,0)
\psline(10,0)(13,0)
\rput(0,-1){
\psdot[linewidth=.05](0,0)
\psdot[linewidth=.05](2.5,0)
\psdot[linewidth=.05](5,0)
\psdot[linewidth=.05](7.5,0)
\psdot[linewidth=.01](9,0)
\psdot[linewidth=.01](9.15,0)
\psdot[linewidth=.01](9.3,0)
\psdot[linewidth=.05](10.5,0)
\psdot[linewidth=.05](13,0)
\rput(0,-.5){
\rput(0,0){$\small{\overline{-g_2(x)}}$}
\rput(2.5,0){$\small{(\overline{g_2(x)}+\overline 1})$}
\rput(5,0){$\small{(-\overline{g_2(x)}+\overline{p-1})}$}
\rput(7.5,0){$\small{(\overline{g_2(x)}+\small{\overline 2})}$}
\rput(10.5,0){$\small{(\overline{g_2(x)}+\overline {\frac{p-1}{2}})}$}
\rput(13,0){$\small{(-\overline{g_2(x)}+\overline{\frac{p+1}{2}})}$}}
\psline(0,0)(8,0)
\psline(10,0)(13,0)}
\psline(0,0)(0,-1)
\psline(13,0)(13,-1)
}
\rput(0,0){
\psdot[linewidth=.05](0,0)
\psdot[linewidth=.05](2.5,0)
\psdot[linewidth=.05](5,0)
\psdot[linewidth=.05](7.5,0)
\psdot[linewidth=.01](9,0)
\psdot[linewidth=.01](9.15,0)
\psdot[linewidth=.01](9.3,0)
\psdot[linewidth=.05](10.5,0)
\psdot[linewidth=.05](13,0)
\rput(0,.3){
\rput(0,0){$\small{\overline{g_{\frac{p^{k-1}-1}{2}}(x)}}$}
\rput(2.5,0){$A$}
\rput(5,0){$B$}
\rput(7.5,0){$C$}
\rput(10.5,0){$D$}
\rput(13,0){$E$}}
\psline(0,0)(8,0)
\psline(10,0)(13,0)

\rput(0,-1){
\psdot[linewidth=.05](0,0)
\psdot[linewidth=.05](2.5,0)
\psdot[linewidth=.05](5,0)
\psdot[linewidth=.05](7.5,0)
\psdot[linewidth=.01](9,0)
\psdot[linewidth=.01](9.15,0)
\psdot[linewidth=.01](9.3,0)
\psdot[linewidth=.05](10.5,0)
\psdot[linewidth=.05](13,0)
\rput(0,-.4){
\rput(0,0){$\small{\overline{-g_{\frac{p^{k-1}-1}{2}}(x)}}$}
\rput(2.5,0){$A^\prime$}
\rput(5,0){$B^\prime$}
\rput(7.5,0){$C^\prime$}
\rput(10.5,0){$D^\prime$}
\rput(13,0){$E^\prime$}}
\psline(0,0)(8,0)
\psline(10,0)(13,0)
\psline(0,0)(0,1)
\psline(13,0)(13,1)}}

\psline[linestyle=dashed,dash=3pt 4pt](3,1.8)(4,1.8)
\psline[linestyle=dashed,dash=3pt 4pt](5,1.8)(6,1.8)
\psline[linestyle=dashed,dash=3pt 4pt](7,1.8)(8,1.8)
\psline[linestyle=dashed,dash=3pt 4pt](3,1.4)(4,1.4)
\psline[linestyle=dashed,dash=3pt 4pt](5,1.4)(6,1.4)
\psline[linestyle=dashed,dash=3pt 4pt](7,1.4)(8,1.4)
\psline[linestyle=dashed,dash=3pt 4pt](3,1)(4,1)
\psline[linestyle=dashed,dash=3pt 4pt](5,1)(6,1)
\psline[linestyle=dashed,dash=3pt 4pt](7,1)(8,1)
}}\end{pspicture}
\caption{Nil clean graph of $GF(p^k)$}\label{gfpk}
\end{figure}

\noindent From the graph, it is clear that the grith of $G_N(\mathbb{F}_p)$ is infinite, which proves $(b)$. It is easy to observe from characterization of finite field that the nil clean graph of $GF(p^k)$ for $ p > 2$, is a disconnected graph consisting of a path of length $p$ and $(\frac{p^{k-1}-1}{2})$ number of $2p$ cycles. For the proof let $GF(p^k) = Z_p[X]/\langle f(x)\rangle$, where $f(x)$ is a irreducible polynomial of degree $k$ over $\mathbb{Z}_{p}$. Let $A\subseteq GF(p^k)$, such that $A$ consists of all linear combinations of $x, x^2, \dots, X^{k-1}$ with coefficients from $\mathbb{Z}_p$ such that if $g(x) + \langle f(x)\rangle \in A$ then $-g(x) + \langle f(x)\rangle \notin A$. Clearly, $A$ can be written as $A=\{g_i(x) + \langle f(x)\rangle ~\big|~ 1\leq i\leq (\frac{p^{k-1}-1}{2})\}$, let us denote $\overline{g_i(x)} = g_i(x) + \langle f(x)\rangle$ for $1 \leq i\leq (\frac{p^{k-1}-1}{2})$. So $(a)$ follows from Figure \ref{gfpk}.

\noindent Here, $A=  -\overline{g_{\frac{p^{k-1}-1}{2}}(x)}+\overline{1}$,
$B=  -\overline{g_{\frac{p^{k-1}-1}{2}}(x)}+\overline{p-1} $,
$C=  \overline{-g_{\frac{p^{k-1}-1}{2}}(x)}+ \overline {2}$,\\

$D=  \overline{-g_{\frac{p^{k-1}-1}{2}}(x)}+\overline {\frac{p-1}{2}} $,
$ E=  -\overline{g_{\frac{p^{k-1}-1}{2}}(x)}+\overline{\frac{p+1}{2}}$,\\

$A^\prime=  \overline{g_{\frac{p^{k-1}-1}{2}}(x)}+\bar{1}$,
$B^\prime=-\overline{g_{\frac{p^{k-1}-1}{2}}(x)}+\overline{p-1} $,
$C^\prime=  \overline{g_{\frac{p^{k-1}-1}{2}}(x)}+ \overline{ 2}$,\\

$D^\prime=  \overline{g_{\frac{p^{k-1}-1}{2}}(x)}+\overline {\frac{p-1}{2}} $, $E^\prime=  -\overline{g_{\frac{p^{k-1}-1}{2}}(x)}+\overline{\frac{p+1}{2}} $

\begin{Cor}
$G_N(R)$ is not cyclic.
\end{Cor}

A  graph $G$ is said to be bipartite if its vertex set can be partitioned into two
disjoint subsets $V_1$ and $V_2$, such that $V(G) = V_1\cup V_2$ and every edge in $G$ has the form $e=(x,y) \in E(G)$, where $x\in V_1$ and $ y \in V_2$. Note that no two vertices both in $V_1$ or both in $V_2$ are adjacent.

\begin{Thm}
$G_N(R)$ is bipartite if and only if $R$ is a field.
\end{Thm}
\noindent$Proof.$ Let $G_N(R)$ be bipartite, so girth is of $G_N(R)$ is not a odd number, hence by  \textbf{Theorem \ref{girth}} $R$ can not be a non field ring. Now if $R$ is a field, it is clear from the nil clean graph of $R$ that $G_N(R)$ is bipartite.$\hfill \square$

\subsection{Dominating set}
Let G be a graph, a subset $S\subseteq V(G)$ is said to be dominating set for $G$ if for all $x\in V(G)$, $x\in S$ or there exists $y\in S$ such that $x$ is adjacent to $y$. Following theorem shows that for a finite commutative weak clean ring dominating number is $2$, where dominating number is the carnality of smallest dominating set.

\begin{Thm}\label{dom_1}
Let $R$ be a weak nil clean ring such that $R$ has no non trivial idempotents, then
$\{ 1, 2\}$ is a dominating set for $G_N(R)$.
\end{Thm}

\noindent$Proof.$ Let $a \in R$, then $a = n, n+1$ or $n-1$ for some $n\in\nil(R)$.
If $a = n,$ then $n+1 \in$ NC($R$) implies $a$ is adjacent to $1$.
If $a = n - 1,$ then $n-1 + 1 = n \in$ NC($R$) implies $a$ is adjacent to $1$.
If $a = n + 1$ and $2 = n_1$ for some nilpotent $n_1 \in R$, then $a$ is adjacent to $2$.
If $a = n + 1$ and $2 = n_1-1$ for some nilpotent $n_1 \in R$, then $a$ is adjacent to $2$.
If $a = n + 1$ and $2 = n_1 + 1$ for some nilpotent $n_1 \in R$. Now
$a + 2 = (n+1) +(n_1+1) = (n+n_1)+2=n+ n_1+(n_1+1) = (n+2n_1) + 1$
is nil clean. Hence $a$ is adjacent to $2$.$\hfill \square$

 \begin{Thm}\label{dom2}
 Let $R = A \times B$, such that $A$ is nil clean and $B$ weak nil clean with no non trivial idempotents. Then $\{(1_A, 1_B), (2_A, 2_B)\}$ is a dominating set for $G_N(R)$.
 \end{Thm}

\noindent$Proof.$ Let $(a,b) \in R$, where $a \in A$ and $b\in B$. For $n_1 \in \nil(A)$, $n_2 \in \nil(B)$ and $0 \neq e\in \idem(A)$, $(a, b)$ has one of the following form
$(a,b) = (n_1, n_2) + (e,1_B), (n_1, n_2) + (e, 0), (n_1, n_2) - (e,1_B)$ or $(n_1, n_2) + (0,0)$
If $(a,b) = (n_1, n_2) + (e,1_B)$, we have $(a, b) + (2_A, 2_B) = (n_1+e+2_A, n_2+1_B+2_B)$. Since $A$ is nil clean $n_1+e+2_A = n^\prime _1 +f$ for some $n_1^\prime \in \nil(A)$ and $f \in \idem(A)$. Since $B$ is weak nil clean $ 2_B =n_2^\prime,$ or $ n_2^\prime - 1_B$, for some $n_2^\prime \in \nil(B)$. If $2_B = n_2^\prime$, we have $(a, b) + (2_A, 2_B) = (n^\prime_1, n_2 + n_2^\prime) + (f, 1_B)$ which is a nil clean expression, hence $(a, b)$ is adjacent to $(2_A, 2_B)$. If $2_B =n_2^\prime - 1_B$, we have
$(a, b) + (2_A, 2_B) = (n^\prime_1, n_2+n^\prime_2) + (f, 0)$, thus $(a, b)$ is adjacent to $(2_A, 2_B)$. In other three cases it is easy to see that $(a,b) + (1_A, 1_B)$ is nil clean hence $(a, b)$ is adjacent to
$(1_A, 1_B)$. Therefore $\{(1_A, 1_B), (2_A, 2_B)\}$ is a dominating set for $R$.$\hfill \square$

\begin{Thm}\label{dom3}
Let $R$ be a weak nil clean ring. $\{ 1, 2\}$ is a dominating set in $G_N(R)$.
\end{Thm}
\noindent$Proof$ If $R$ has no non trivial idempotents, then by \textbf{Theorem \ref{dom_1}} we are done.
If $R$ has a non trivial idempotent say $e$, then by Peirce decomposition
$R\cong eR \oplus (1-e)R$, now by \textbf{Theorem 2.3} of \cite{wnc} we have one of $eR$ or $(1-e)R$ must be a nil clean ring. Without lost of generality suppose that $eR$ is a nil clean ring and $(1-e)R$ be a weak nil clean ring. Now if $(1-e)R$ has no non trivial idempotents, then we have the result by \textbf{Theorem \ref{dom2}}. If $f \in \idem((1-e)R)$, repeating as above we get a direct sum decomposition of $R$ where only one summand is weak nil clean. As ring $R$ is a finite ring, so after repeating above to the weak nil clean summand of $R$, we will have a direct sum decomposition of $R$, where idempotents of weak nil clean summand of $R$ is trivial, then again by \textbf{Theorem \ref{dom2}} we have the result.$\hfill \square$
\subsection{Chromatic index}
 An edge colouring of a graph $G$  is a map  $C: E(G) \rightarrow S$, where $S$ is a set of colours such that for all $e, f \in E(G)$, if $e$ and $f$ are adjacent, then $C(e) \neq C(f)$. The \textit{chromatic index} of a graph denoted by $\chi^\prime(G)$ and is defined as the minimum number of colours needed for a proper colouring of $G$. Let $\vartriangle$ be the maximum vertex degree of $G$, then Vizing’s theorem \cite{gt} gives $\vartriangle \leq \chi^\prime(G) \leq \vartriangle+ 1$. Vizing’s theorem divides the graphs into two classes according to their chromatic index; graphs satisfying $\chi^\prime(G) = \vartriangle$ are called graphs of class $1$, those with $\chi^\prime(G) =\vartriangle + 1$ are graphs of class $2$.

Following theorem shows that $G_N(R)$ is of class $1$.

\begin{Thm}
Let $R$ be a finite commutative ring then the nil clean graph of $R$ is of class $1$.
\end{Thm}

\noindent$Proof.$ We colour the edge $ab$ by the colour $a + b$. By this colouring, every two distinct edges $ab$ and $ac$ has a different colour and $C = \{ a + b| ab\mbox{ is an edge in }G_N(R) \}$ is the set of colours. Therefore nil clean graph has a $|C|$-edge colouring and so $\chi^\prime(G_N(R)) \leq |C|$. But
$C \subset NC(R)$ and $\chi^\prime(G_N(R)) \leq |C| \leq |NC(R)|$. By \textbf{Lemma \ref{deg}},  $\vartriangle \leq |NC(R)|,$ and so by Vizing's theorem, we have $\chi^\prime(G_N(R))\geq \vartriangle = |NC(R|$. Therefore $\chi^\prime(G_N(R)) = |NC(R| =\vartriangle,$ i.e., $G_N(R)$ is class of $1$.$\hfill \square$
\subsection{Diameter}
For a graph $G$, number of edges on the shortest path between vertices $x$ and $y$ is called the \textit{distance} between $x$ and $y$ and is denoted by $d(x,y)$. If there is no path between $x$ and $y$ then we say $d(x,y)= \infty$. The diameter of a graph $diam(G)$ is the maximum of distances of each pair of distinct vertices in $G$.

The following are some results related to diameter of nil clean graph of a ring.

\begin{Lem}\label{dia1}
$R$ is nil clean ring if and only if $diam(G_N(R))=1$.
\end{Lem}

\begin{Thm}
Let $R$ be a non nil clean, weak nil clean ring with no non trivial idempotents then $diam(G_N(R)) = 2$.
\end{Thm}

\noindent$Proof.$ Let $a, b \in R$, then for some $n_1, n_2 \in \nil(R)$, we have $a = n_1 + 1, n_1 -1$ or $n_1$ and $b = n_2 + 1, n_2 -1$ or $n_2$. If $a = n_1 + 1$ and $b =  n_2 -1$ or $n_2$, clearly $a+b$ is nil clean, thus $ab$ is an edge in $G_N(R)$, therefore $d(a,b) = 1$. If $a = n_1 + 1$ and $b =  n_2 +1$, we have the path $a \rule[0.5ex]{.5cm}{0.4pt} (-1) \rule[0.5ex]{.5cm}{0.4pt} b $ in $G_N(R)$ thus $d(a,b)\leq 2$. If $a = n_1 -1$ and $b = n_2 - 1$ or $n_2$ similarly as above a path of length $2$ from $a$ to $b$ through $1$, thus in this case also $d(a,b)\leq 2$. Finally if $a = n_1$ and $b=n_2$, then $d(a,b) = 2$.
Thus from above we conclude that $diam(G_N(R))\leq2$. Now as $R$ is a non nil clean, weak nil clean ring thus we have at least one $x\in R$, such that $x = n-1$ but $x\neq n +1$, i.e., $x$ is not nil clean. Then we have $d(0,x) = 2$, therefore $diam(G_N(R)) \geq 2$, and hence the result follows.$\hfill \square$
\begin{Thm}\label{dia2}
Let $R = A \times B$, such that $A$ is nil clean and $B$ weak nil clean with no non trivial idempotents,
then $diam(G_N(R)) = 2$.
\end{Thm}

\noindent$Proof.$  We have $\idem(R) = \{(e,0_B), (e,1_B) | e\in \idem(A)\}$. Now let
$(a_1,b_1), (a_2,b_2) \in R$, in case $(a_1,b_1) + (a_2,b_2)$ is nil clean then
$d((a_1,b_1), (a_2,b_2) ) =  1$ in $G_N(R)$. If $(a_1,b_1) + (a_2,b_2)$ is not nil clean, thus $b_1 +b_2$ is not nil clean. So we have followings cases:\\
\textbf{CASE I:} If $b_1 = n_1 +1$ and $b_2 = n_2 +1$, we have the path $(a_1,b_1)\rule[0.5ex]{.5cm}{0.4pt} (0,-1) \rule[0.5ex]{.5cm}{0.4pt} (a_2,b_2)  $ in $G_N(R)$, thus $d((a_1,b_1),(a_2,b_2))\leq 2$.\\
\textbf{CASE II:} If $b_1 = n_1 -1$ and $b_2 = n_2 -1$, we have the path $(a_1,b_1)\rule[0.5ex]{.5cm}{0.4pt} (0,1) \rule[0.5ex]{.5cm}{0.4pt} (a_2,b_2)  $ in $G_N(R)$, thus $d((a_1,b_1),(a_2,b_2))\leq 2$.\\
\textbf{CASE III:} If $b_1 = n_1 -1$ and $b_2 = n_2$, we have the path $(a_1,b_1)\rule[0.5ex]{.5cm}{0.4pt} (0,1) \rule[0.5ex]{.5cm}{0.4pt} (a_2,b_2)  $ in $G_N(R)$, thus $d((a_1,b_1),(a_2,b_2))\leq 2$.\\
\textbf{CASE IV:} If $b_1 = n_1$ and $b_2 = n_2 -1$, by \textbf{CASE III} and symmetry we have $d((a_1,b_1),(a_2,b_2))\leq 2$.\\
Therefore $diam(R) \leq 2$, $R$ is not nil clean implies $diam(R) \geq 2$, Thus $diam(R) = 2$.$\hfill \square$

\begin{Thm}\label{dia3}
 If $R$ is weak nil clean ring but not nil clean then $diam(G_N(R))=2$.
\end{Thm}

\noindent$Proof.$ If $R$ has no non trivial idempotents, then by \textbf{Theorem \ref{dia2}} we are done.
If $R$ has a non trivial idempotent say $e$, then by Peirce decomposition
$R\cong eR \oplus (1-e)R$, so by \textbf{Theorem 2.3} of \cite{wnc} we have one of $eR$ or $(1-e)R$ must be a nil clean ring. Without lost of generality suppose that $eR$ is a nil clean ring and $(1-e)R$ be a weak nil clean ring. Now if $(1-e)R$ has no non trivial idempotents, then we have the result by \textbf{Theorem \ref{dia2}}. If $f \in \idem((1-e)R)$, repeating as above we get a direct sum decomposition of $R$ where only one summand is weak nil clean. As $R$ is a finite ring, so after repeating above to the weak nil clean summand of $R$ we will have a direct sum decomposition of $R$, where idempotents of weak nil clean summand of $R$ are trivial, then again by \textbf{Theorem \ref{dia2}} we have the result.$\hfill \square$

\begin{Thm}
Let $n$ be a positive integer, then the following hold for $\mathbb{Z}_n$.
\begin{enumerate}
  \item If $n = 2^k$, for some integer $k\geq 1$, then $diam(G_N(\mathbb{Z}_n) )= 1$.
  \item If $n = 2^k3^l$, for some integer $k\geq 0$ and $l \geq 1$, then $diam(G_N(\mathbb{Z}_n)) = 2$.
  \item  For a prime $p$, $diam(G_N(\mathbb{Z}_p)) = p-1$.
  \item If $n=2p$, where $p$ is a odd prime, then $diam(G_N(\mathbb{Z}_{2p}) )= p-1$.
  \item If $n=3p$, where $p$ is a odd prime, then $diam(G_N(\mathbb{Z}_{3p})) = p-1$.
\end{enumerate}
\end{Thm}

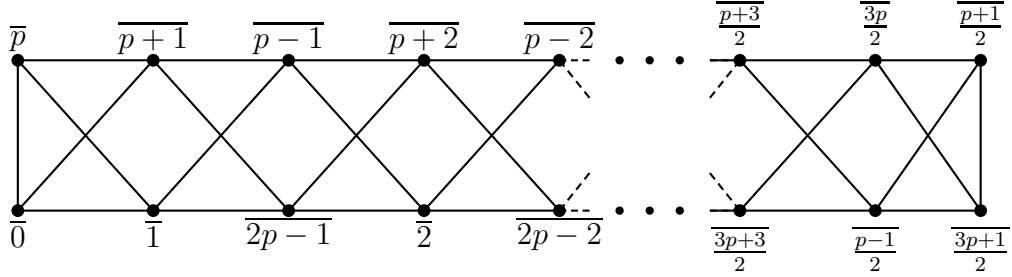
\begin{figure}[H]   
\centering
\begin{pspicture}(0,2)(0,-2)
\scalebox{1}{
\rput(-6.5,0){
\psdot[linewidth=.05](0,1.5)
\psdot[linewidth=.05](1.8,1.5)
\psdot[linewidth=.05](3.6,1.5)
\psdot[linewidth=.05](5.4,1.5)
\psdot[linewidth=.05](7.2,1.5)
\psdot[linewidth=.02](8,1.5)
\psdot[linewidth=.02](8.4,1.5)
\psdot[linewidth=.02](8.8,1.5)
\psdot[linewidth=.05](9.6,1.5)
\psdot[linewidth=.05](11.4,1.5)
\psdot[linewidth=.05](12.8,1.5)
\rput(0,1.8){$\overline {p}$}
\rput(1.8,1.8){$\overline {p+1}$}
\rput(3.6,1.8){$\overline {p-1}$}
\rput(5.4,1.8){$\overline {p+2}$}
\rput(7.2,1.8){$\overline {p-2}$}
\rput(9.6,2){$\overline {\frac{p+3}{2}}$}
\rput(11.4,2){$\overline {\frac{3p}{2}}$}
\rput(12.8,2){$\overline {\frac{p+1}{2}}$}
\psdot[linewidth=.05](0,-.5)
\psdot[linewidth=.05](1.8,-.5)
\psdot[linewidth=.05](3.6,-.5)
\psdot[linewidth=.05](5.4,-.5)
\psdot[linewidth=.05](7.2,-.5)
\psdot[linewidth=.02](8,-.5)
\psdot[linewidth=.02](8.4,-.5)
\psdot[linewidth=.02](8.8,-.5)
\psdot[linewidth=.05](9.6,-.5)
\psdot[linewidth=.05](11.4,-.5)
\psdot[linewidth=.05](12.8,-.5)
\rput(0,-.8){$\overline {0}$}
\rput(1.8,-.8){$\overline {1}$}
\rput(3.6,-.8){$\overline {2p-1}$}
\rput(5.4,-.8){$\overline {2}$}
\rput(7.2,-.8){$\overline {2p-2}$}
\rput(9.6,-1){$\overline {\frac{3p+3}{2}}$}
\rput(11.4,-1){$\overline {\frac{p-1}{2}}$}
\rput(12.8,-1){$\overline {\frac{3p+1}{2}}$}
\psline(0,1.5)(7.2,1.5)
\psline(0,1.5)(0,-.5)
\psline(0,-.5)(7.2,-.5)
\psline(9.2,1.5)(12.8,1.5)(12.8,-.5)(9.2,-.5)
\psline(0,1.5)(1.8,-.5)(3.6,1.5)(5.4,-.5)(7.2,1.5)
\psline(0,-.5)(1.8,1.5)(3.6,-.5)(5.4,1.5)(7.2,-.5)
\psline[linestyle=dashed,dash=3pt 2pt](7.2,1.5)(7.6,1.5)
\psline[linestyle=dashed,dash=3pt 2pt](7.2,1.5)(7.6,1)
\psline[linestyle=dashed,dash=3pt 2pt](7.2,-.5)(7.6,0)
\psline[linestyle=dashed,dash=3pt 2pt](7.2,-.5)(7.6,-.5)
\psline(9.6,1.5)(11.4,-.5)(12.8,1.5)
\psline(9.6,-.5)(11.4,1.5)(12.8,-.5)
\psline[linestyle=dashed,dash=3pt 2pt](9.2,1)(9.6,1.5)
\psline[linestyle=dashed,dash=3pt 2pt](9.2,0)(9.6,-.5)
\psline[linestyle=dashed,dash=3pt 2pt](9.2,1.5)(9.6,1.5)
\psline[linestyle=dashed,dash=3pt 2pt](9.2,-.5)(9.6,-.5)
}}
\end{pspicture}
\caption{Nil clean graph of $\mathbb{Z}_{2p}$ }\label{2p}
\end{figure}

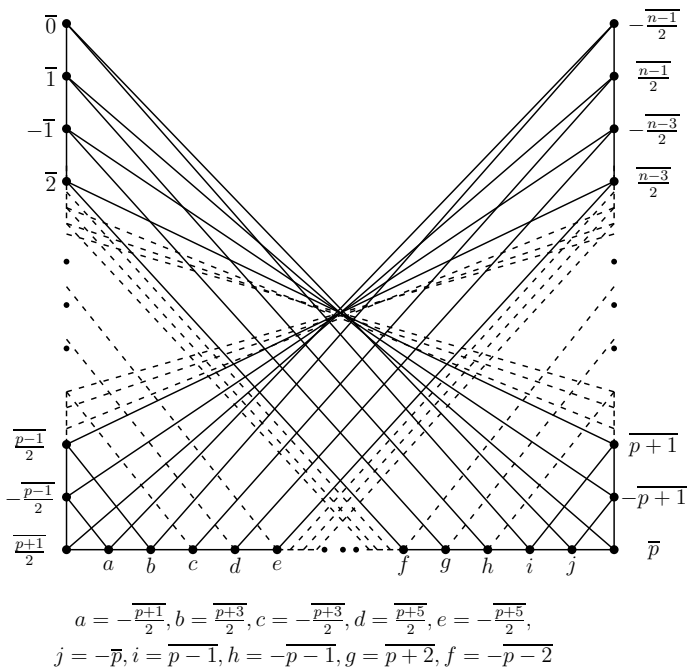
\begin{figure}[H]
\begin{pspicture}(0,-2.3)(0,7.2)
\scalebox{.7}{

\rput(3.5,0){

\psdot[linewidth=.05](0,0)
\psdot[linewidth=.05](.8,0)
\psdot[linewidth=.05](1.6,0)
\psdot[linewidth=.05](2.4,0)
\psdot[linewidth=.05](3.2,0)
\psdot[linewidth=.05](4,0)
\psdot[linewidth=.02](4.9,0)
\psdot[linewidth=.02](5.25,0)
\psdot[linewidth=.02](5.5,0)
\psdot[linewidth=.05](6.4,0)
\psdot[linewidth=.05](7.2,0)
\psdot[linewidth=.05](8,0)
\psdot[linewidth=.05](8.8,0)
\psdot[linewidth=.05](9.6,0)
\psdot[linewidth=.05](10.4,0)

\psdot[linewidth=.05](10.4,1)
\psdot[linewidth=.05](10.4,2)
\psdot[linewidth=.02](10.4,3.825)
\psdot[linewidth=.02](10.4,4.65)
\psdot[linewidth=.02](10.4,5.475)

\psdot[linewidth=.05](10.4,7)
\psdot[linewidth=.05](10.4,8)
\psdot[linewidth=.05](10.4,9)
\psdot[linewidth=.05](10.4,10)

\psdot[linewidth=.05](0,1)
\psdot[linewidth=.05](0,2)
\psdot[linewidth=.02](0,3.825)
\psdot[linewidth=.02](0,4.65)
\psdot[linewidth=.02](0,5.475)
\psdot[linewidth=.05](0,7)
\psdot[linewidth=.05](0,8)
\psdot[linewidth=.05](0,9)
\psdot[linewidth=.05](0,10)
\psline(.8,0)(10.4,10)
\psline(1.6,0)(10.4,10)
\psline(0,0)(10.4,9)
\psline(2.4,0)(10.4,9)
\psline(0,1)(10.4,8)
\psline(3.2,0)(10.4,8)
\psline(0,2)(10.4,7)
\psline(4,0)(10.4,7)

\psline[linestyle=dashed,dash=3pt 4pt](0,3)(10.4,6)
\psline[linestyle=dashed,dash=3pt 4pt](0,2.7)(10.4,6.2)
\psline[linestyle=dashed,dash=3pt 4pt](0,2.3)(10.4,6.5)

\psline[linestyle=dashed,dash=3pt 4pt](10.4,3)(0,6)
\psline[linestyle=dashed,dash=3pt 4pt](10.4,2.7)(0,6.2)
\psline[linestyle=dashed,dash=3pt 4pt](10.4,2.3)(0,6.5)

\psline[linestyle=dashed,dash=3pt 4pt](4.25,0)(10.4,6.8)
\psline[linestyle=dashed,dash=3pt 4pt](4.5,0)(10.4,6.5)
\psline[linestyle=dashed,dash=3pt 4pt](4.75,0)(10.4,6.2)

\psline[linestyle=dashed,dash=3pt 4pt](6.15,0)(0,6.8)
\psline[linestyle=dashed,dash=3pt 4pt](5.9,0)(0,6.5)
\psline[linestyle=dashed,dash=3pt 4pt](5.75,0)(0,6.2)

\psline(8.8,0)(0,10)
\psline(9.6,0)(0,10)
\psline(10.4,0)(0,9)
\psline(8,0)(0,9)
\psline(10.4,1)(0,8)
\psline(7.2,0)(0,8)
\psline(10.4,2)(0,7)
\psline(6.4,0)(0,7)

\psline(0,1)(.8,0)
\psline(0,2)(1.6,0)
\psline[linestyle=dashed,dash=3pt 4pt](0,3)(2.4,0)
\psline[linestyle=dashed,dash=3pt 4pt](0,4)(3.2,0)
\psline[linestyle=dashed,dash=3pt 4pt](0,5)(4,0)

\psline(10.4,1)(9.6,0)
\psline(10.4,2)(8.8,0)
\psline[linestyle=dashed,dash=3pt 4pt](10.4,3)(8,0)
\psline[linestyle=dashed,dash=3pt 4pt](10.4,4)(7.2,0)
\psline[linestyle=dashed,dash=3pt 4pt](10.4,5)(6.4,0)

\psline[linestyle=dashed,dash=3pt 2pt](0,2)(0,3)
\psline(0,2)(0,0)(4,0)
\psline[linestyle=dashed,dash=3pt 2pt](4,0)(4.7,0)
\psline[linestyle=dashed,dash=3pt 2pt](5.7,0)(6.4,0)
\psline(6.4,0)(10.4,0)(10.4,2)
\psline[linestyle=dashed,dash=3pt 2pt](10.4,2)(10.4,3)

\psline(10.4,7)(10.4,10)
\psline[linestyle=dashed,dash=3pt 2pt](10.4,6.3)(10.4,7.3)

\psline(0,7)(0,10)
\psline[linestyle=dashed,dash=3pt 2pt](0,6.3)(0,7.3)

\rput(-.7,0){$\small{\overline {\frac{p+1}{2}}}$}
\rput(-.7,1){$\small{-\overline {\frac{p-1}{2}}}$}
\rput(-.7,2){$\small{ \overline {\frac{p-1}{2}}}$}
\rput(-.3,10){$\small{\overline {0}}$}
\rput(-.3,9){$\small{\overline {1}}$}
\rput(-.5,8){$\small{-\overline {1}}$}
\rput(-.3,7){$\small{\overline {2}}$}
\rput(0.8,-.3){$\small{a}$}
\rput(1.6,-.3){$\small{b}$}
\rput(2.4,-.3){$\small{c}$}
\rput(3.2,-.3){$\small{d}$}
\rput(4,-.3){$\small{e}$}
\rput(6.4,-.3){$\small{f}$}
\rput(7.2,-.3){$\small{g}$}
\rput(8,-.3){$\small{h}$}
\rput(8.8,-.3){$\small{i}$}
\rput(9.6,-.3){$\small{j}$}
\rput(4.5,-1.3){$\small{a = -\overline {\frac{p+1}{2}}, b=\overline {\frac{p+3}{2}},
c =-\overline {\frac{p+3}{2}},d =\overline {\frac{p+5}{2}},e =-\overline {\frac{p+5}{2}},}$}
 \rput(4.5,-2){$\small{j =-\overline {p},i =\overline {p-1},h =-\overline {p-1},g=\overline {p+2},
 f =-\overline {p-2}}$}

 \rput(11.15,0){$\small{\overline {p}}$}
\rput(11.15,1){$\small{-\overline {p+1}}$}
\rput(11.15,2){$\small{\overline {p+1}}$}
\rput(11.15,10){$\small{-\overline {\frac{n-1}{2}}}$}
\rput(11.15,9){$\small{\overline {\frac{n-1}{2}}}$}
\rput(11.15,8){$\small{-\overline {\frac{n-3}{2}}}$}
\rput(11.15,7){$\small{\overline {\frac{n-3}{2}}}$}
}}
\end{pspicture}
\caption{Nil clean graph of $\mathbb{Z}_{3p}$ where $p \equiv 2(mod 3)$}\label{3p2}
\end{figure}

\noindent$Proof.$ $(i)$ and $(ii)$ follow from \textbf{Lemma \ref{dia1}} and \textbf{Theorem \ref{dia3}} respectively. $(iii)$ follows from $G_N(\mathbb{Z}_{p})$ in Figure \ref{zp}, $(iv)$  follows from graph  $G_N(\mathbb{Z}_{2p})$  in Figure \ref{2p} and $(v)$ follows from graphs in Figure \ref{3p2} and Figure \ref{3p1}.

\begin{figure}[H]
\begin{pspicture}(0,-2.3)(0,7.2)
\scalebox{.7}{

\rput(3.5,0){
\psdot[linewidth=.05](0,0)
\psdot[linewidth=.05](.8,0)
\psdot[linewidth=.05](1.6,0)
\psdot[linewidth=.05](2.4,0)
\psdot[linewidth=.05](3.2,0)
\psdot[linewidth=.05](4,0)
\psdot[linewidth=.02](4.9,0)
\psdot[linewidth=.02](5.25,0)
\psdot[linewidth=.02](5.5,0)
\psdot[linewidth=.05](6.4,0)
\psdot[linewidth=.05](7.2,0)
\psdot[linewidth=.05](8,0)
\psdot[linewidth=.05](8.8,0)
\psdot[linewidth=.05](9.6,0)
\psdot[linewidth=.05](10.4,0)
\psdot[linewidth=.05](10.4,1)
\psdot[linewidth=.05](10.4,2)
\psdot[linewidth=.02](10.4,3.825)
\psdot[linewidth=.02](10.4,4.65)
\psdot[linewidth=.02](10.4,5.475)
\psdot[linewidth=.05](10.4,7)
\psdot[linewidth=.05](10.4,8)
\psdot[linewidth=.05](10.4,9)
\psdot[linewidth=.05](10.4,10)

\psdot[linewidth=.05](0,1)
\psdot[linewidth=.05](0,2)
\psdot[linewidth=.02](0,3.825)
\psdot[linewidth=.02](0,4.65)
\psdot[linewidth=.02](0,5.475)
\psdot[linewidth=.05](0,7)
\psdot[linewidth=.05](0,8)
\psdot[linewidth=.05](0,9)
\psdot[linewidth=.05](0,10)
\psline(.8,0)(10.4,10)
\psline(1.6,0)(10.4,10)
\psline(0,0)(10.4,9)
\psline(2.4,0)(10.4,9)
\psline(0,1)(10.4,8)
\psline(3.2,0)(10.4,8)
\psline(0,2)(10.4,7)
\psline(4,0)(10.4,7)

\psline[linestyle=dashed,dash=3pt 4pt](0,3)(10.4,6)
\psline[linestyle=dashed,dash=3pt 4pt](0,2.7)(10.4,6.2)
\psline[linestyle=dashed,dash=3pt 4pt](0,2.3)(10.4,6.5)

\psline[linestyle=dashed,dash=3pt 4pt](10.4,3)(0,6)
\psline[linestyle=dashed,dash=3pt 4pt](10.4,2.7)(0,6.2)
\psline[linestyle=dashed,dash=3pt 4pt](10.4,2.3)(0,6.5)

\psline[linestyle=dashed,dash=3pt 4pt](4.25,0)(10.4,6.8)
\psline[linestyle=dashed,dash=3pt 4pt](4.5,0)(10.4,6.5)
\psline[linestyle=dashed,dash=3pt 4pt](4.75,0)(10.4,6.2)

\psline[linestyle=dashed,dash=3pt 4pt](6.15,0)(0,6.8)
\psline[linestyle=dashed,dash=3pt 4pt](5.9,0)(0,6.5)
\psline[linestyle=dashed,dash=3pt 4pt](5.75,0)(0,6.2)

\psline(8.8,0)(0,10)
\psline(9.6,0)(0,10)
\psline(10.4,0)(0,9)
\psline(8,0)(0,9)
\psline(10.4,1)(0,8)
\psline(7.2,0)(0,8)
\psline(10.4,2)(0,7)
\psline(6.4,0)(0,7)

\psline(0,1)(.8,0)
\psline(0,2)(1.6,0)
\psline[linestyle=dashed,dash=3pt 4pt](0,3)(2.4,0)
\psline[linestyle=dashed,dash=3pt 4pt](0,4)(3.2,0)
\psline[linestyle=dashed,dash=3pt 4pt](0,5)(4,0)

\psline(10.4,1)(9.6,0)
\psline(10.4,2)(8.8,0)
\psline[linestyle=dashed,dash=3pt 4pt](10.4,3)(8,0)
\psline[linestyle=dashed,dash=3pt 4pt](10.4,4)(7.2,0)
\psline[linestyle=dashed,dash=3pt 4pt](10.4,5)(6.4,0)

\psline[linestyle=dashed,dash=3pt 2pt](0,2)(0,3)
\psline(0,2)(0,0)(4,0)
\psline[linestyle=dashed,dash=3pt 2pt](4,0)(4.7,0)
\psline[linestyle=dashed,dash=3pt 2pt](5.7,0)(6.4,0)
\psline(6.4,0)(10.4,0)(10.4,2)
\psline[linestyle=dashed,dash=3pt 2pt](10.4,2)(10.4,3)

\psline(10.4,7)(10.4,10)
\psline[linestyle=dashed,dash=3pt 2pt](10.4,6.3)(10.4,7.3)

\psline(0,7)(0,10)
\psline[linestyle=dashed,dash=3pt 2pt](0,6.3)(0,7.3)

\rput(-.7,0){$\small{- \overline {\frac{p-1}{2}}}$}
\rput(-.7,1){$\small{\overline {\frac{p-1}{2}}}$}
\rput(-.7,2){$\small{- \overline {\frac{p-3}{2}}}$}
\rput(-.3,10){$\small{\overline {0}}$}
\rput(-.3,9){$\small{\overline {1}}$}
\rput(-.5,8){$\small{-\overline {1}}$}
\rput(-.3,7){$\small{\overline {2}}$}
\rput(0.8,-.3){$\small{a}$}
\rput(1.6,-.3){$\small{b}$}
\rput(2.4,-.3){$\small{c}$}
\rput(3.2,-.3){$\small{d}$}
\rput(4,-.3){$\small{e}$}
\rput(6.4,-.3){$\small{f}$}
\rput(7.2,-.3){$\small{g}$}
\rput(8,-.3){$\small{h}$}
\rput(8.8,-.3){$\small{i}$}
\rput(9.6,-.3){$\small{j}$}
\rput(4.5,-1.3){$\small{a = \overline {\frac{p+1}{2}}, b=-\overline {\frac{p+1}{2}},
c =\overline {\frac{p+3}{2}},d =-\overline {\frac{p+3}{2}},e =\overline {\frac{p+5}{2}},}$}
 \rput(4.5,-2){$\small{j =\overline {p},i =-\overline {p-1},h =\overline {p-1},g=-\overline {p+2},
 f =\overline {p-2}}$}

 \rput(11.15,0){$\small{- \overline {p}}$}
\rput(11.15,1){$\small{\overline {p+1}}$}
\rput(11.15,2){$\small{-\overline {p+1}}$}
\rput(11.15,10){$\small{-\overline {\frac{n-1}{2}}}$}
\rput(11.15,9){$\small{\overline {\frac{n-1}{2}}}$}
\rput(11.15,8){$\small{-\overline {\frac{n-3}{2}}}$}
\rput(11.15,7){$\small{\overline {\frac{n-3}{2}}}$}

}}
\end{pspicture}
\caption{Nil clean graph of $\mathbb{Z}_{3p}$ where $p\equiv (1mod 3)$}\label{3p1}
\end{figure}


\end{document}